\documentclass{french_current}
\usepackage{latexsym}
\usepackage[french]{babel}
\RequirePackage{nccmath}
\RequirePackage[colorlinks, bookmarks=true]{hyperref}

\newcommand{\noun}[1]{\textsc{#1}}

\begin{document}
\subh{Smooths Tests of Goodness-of-fit for the Newcomb-Benford distribution }
\author{Gilles R. DUCHARME\ts{1},Samuel KACI\ts{1},
Credo VOVOR-DASSU\ts{1}}

\aff{\removelastskip \ts{1}IMAG, Univ. Montpellier, CNRS, Montpellier, France}

\resu{La loi de probabilit\'{e} de Newcomb-Benford est de plus en plus utilis\'{e}e
dans les applications de la statistique, notamment en d\'{e}tection de fraude. Dans ces contextes, il importe de d\'{e}terminer si un jeu de donn\'{e}es est issu de cette loi de probabilit\'{e} en contrôlant les risques d'erreur de Type 1, soit de faussement identifier une fraude, et de Type 2, soit de ne pas la d\'{e}tecter. L'outil statistique qui permet d'ex\'{e}cuter ce genre de t\^{a}che est le test d'ad\'{e}quation. Pour la loi de Newcomb-Benford, le test d'ad\'{e}quation le plus populaire est le test du khi-deux de Pearson dont la puissance, associ\'{e}e au risque d'erreur de Type 2, est reconnue comme \'{e}tant assez faible. En cons\'{e}quence, d'autres tests ont \'{e}t\'{e} r\'{e}cemment introduits. Le but de ce travail est de proposer de nouveaux tests d'ad\'{e}quation pour cette loi, bas\'{e}s sur le principe des tests lisses. Ces tests sont ensuite compar\'{e}s aux meilleurs tests existants pour ce probl\`{e}me. Il en ressort que nos propositions sont globalement pr\'{e}f\'{e}rables aux tests existants et pourraient \^{e}tre utilis\'{e}es dans les applications, notamment en d\'{e}tection de fraude. Un package de R \noun{BenfordSmoothTest} est disponible sur le site \href{https://github.com/st-homme/benfordSmoothTest}{GitHub} pour calculer les statistiques de test.}

\abstract{The Newcomb-Benford probability distribution is becoming very popular
in many areas using statistics, notably in fraud detection. In such contexts, it is important to be able to determine if a data set arises from this distribution while controlling the risk of a Type 1 error, i.e. falsely identifying a fraud, and a Type 2 error, i.e. not detecting that a fraud occurred. The statistical tool to do this work is a goodness-of-fit test. For the Newcomb-Benford distribution, the most popular such test is Pearson's chi-square test whose power, related to the Type
2 error, is known to be weak. Consequently, other tests have been recently introduced. The goal of the present work is to build new goodness-of-fit tests for this distribution, based on the smooth test principle. These tests are then compared to some of their competitors. It turns out that the proposals of the paper are globally preferable to existing tests and should be seriously considered in fraud detection contexts, among others. The R package \noun{BenfordSmoothTest} is available on \href{https://github.com/st-homme/benfordSmoothTest}{GitHub} to compute the test statistics}

\moi{\textbf{\textcolor{AbsBlue}{MOTS-CL\'{E}S.}}\hphantom{--} Loi de Newcomb-Benford,  test d'ad\'{e}quation, d\'{e}tection de fraude, test lisse}
\kwd{Newcomb-Benford's distribution, goodness-of-fit tests, faud detection; smooth test}

\chapter{Tests d'ad\'{e}quations lisses pour la loi de Newcomb-Benford}

\section{Introduction \vspace*{-5truemm}}
\rmfamily

La loi de Newcomb-Benford LNB (NEWCOMB 1881), (BENFORD 1938) annonce que sous certaines conditions, le premier chiffre significatif (\emph{PCS}) d'une variable al\'{e}atoire continue positive $X$, $D=PCS(X)$, a pour probabilit\'{e} $\mathbb{P}[D=d]=\pi_{d}=\log_{10}[1+1/d]$, $d\in\{1,2,...,9\}$.

L'utilisation de cette loi de probabilit\'{e} connaît une popularit\'{e} grandissante dans de nombreux domaines, notamment en d\'{e}tection de fraudes fiscales (NIGRINI 1996), financi\`{e}res (CERIOLI et collab. 2019), comptables (DURTSCHI et collab. 2004) et scientifiques (HEIN et collab. 2012). 
Elle est aussi utilis\'{e}e comme mod\`{e}le statistique dans des disciplines
aussi vari\'{e}es que l'hydrologie (NIGRINI et DRAKE 2007), la volcanologie (GEYER et MARTI 2012), la sismologie (SAMBRIDGE et collab. 2011)  et pour l'\'{e}tude du trafic de donn\'{e}es internet (ARSHADI et JAHANGIR 2014), entre autres applications. Cette popularit\'{e} \'{e}mane d'une part, du fait qu'on la rencontre empiriquement tr\`{e}s souvent dans les jeux de donn\'{e}es r\'{e}elles (DURTSCHI et collab. 2004) correspondant \`{a} des donn\'{e}es dites de \og deuxi\`{e}me g\'{e}n\'{e}ration \fg{ }, soit le r\'{e}sultat d'op\'{e}rations (produits, puissances, etc.) de donn\'{e}es brutes. D'autre part, des raisons th\'{e}oriques font qu'elle appara\^{i}t aussi dans de nombreux contextes (POSCH 2008), sinon exactement du moins en tant qu'approximation de la r\'{e}alit\'{e}. En outre, elle interpelle l'intuition parce que, contrairement \`{a} ce que l'on pourrait penser, $D$ n'appara\^{i}t pas avec une probabilit\'{e} de 1/9 $\simeq0.111$. Cette dissonance cognitive vient de ce que les psychologues appellent le \og biais d'\'{e}quiprobabilit\'{e} \fg{ } (LECOUTRE 1992). Ce biais est un des facteurs faisant que les \emph{PCS} de nombres influenc\'{e}s par la pens\'{e}e humaine sont plus proches de la loi uniforme sur $\{1,2,...,9\}$ que de la LNB (HILL 1998), (GAUVRIT et collab. 2017). En particulier un fraudeur voulant trafiquer un jeu de donn\'{e}es va inconsciemment avoir tendance \`{a} uniformiser leur \emph{PCS}.
Cette discordance offre une prise permettant \`{a} des auditeurs de d\'{e}tecter
les fraudes, et cet outil s'est av\'{e}r\'{e} particuli\`{e}rement efficace dans
le cas de fraudes fiscales (AUSLOOS et collab. 2017). Un site internet, le \emph{Benford online bibliograph}y (BERGER et collab. 2015) recense la quasi-totalit\'{e} des publications sur la LNB, autant les r\'{e}sultats th\'{e}oriques que les applications.

Il est donc souvent imp\'{e}ratif de pouvoir d\'{e}terminer si un jeu de donn\'{e}es se conforme \`{a} la LNB en contrôlant les risques d'erreur de Type I
et II. Dans le contexte de la d\'{e}tection de fraudes, ces erreurs correspondent
\`{a} faussement suspecter une fraude (Type I) avec pour cons\'{e}quence le
coût d'une audition approfondie subs\'{e}quente, soit de ne pas d\'{e}tecter
un jeu de donn\'{e}es trafiqu\'{e}es (Type II) qui m\`{e}nera ensuite \`{a} la prise de d\'{e}cisions erron\'{e}es. L'outil statistique pour effectuer cette t\^{a}che est un test d'ad\'{e}quation statistique (\emph{goodness-of-fit test}). 
Pour ce probl\`{e}me, le test d'ad\'{e}quation de pr\'{e}f\'{e}rence  (MORROW 2014) a longtemps \'{e}t\'{e} le test du $\chi^{2}$ de Pearson (LESPERANCE et collab.  2016, eq. 3). Si $n_{d}$ est le nombre de fois dans un jeu de $n$ donn\'{e}es que \emph{$PCS=d$}, alors le test s'effectue en calculant la statistique de test :
\begin{ceqn}
\begin{align}
\chi^{2} & =\sum_{i=1}^{9}\frac{\left(n_{d}-n\pi_{d}\right)^{2}}{n\pi_{d}},\label{eq:Pearson chi-square test}
\end{align}
\end{ceqn}
qui ob\'{e}it approximativement, si $n$ est grand, \`{a} la loi $\chi_{8}^{2}$
si l'hypoth\`{e}se nulle $H_{0}$ voulant que les donn\'{e}es suivent une
LNB est juste. On rejette $H_{0}$ si $\chi^{2}$ d\'{e}passe le quantile
de cette loi appropri\'{e} au risque d'erreur de Type I souhait\'{e}. Ce test
est simple d'utilisation, mais sa puissance n'est en g\'{e}n\'{e}ral pas reconnue comme \'{e}tant tr\`{e}s
\'{e}lev\'{e}e (MORROW 2014). Ainsi d'autres tests d'ad\'{e}quation ont \'{e}t\'{e} propos\'{e}s, certains \'{e}tant des adaptations de tests d\'{e}velopp\'{e}s pour des donn\'{e}es continues et bas\'{e}s sur des principes statistiques reconnus (test de Carmer-von Mises ou de Watson dans LESPERANCE et collab. 2016), d'autres sur la base de consid\'{e}rations intuitives (JOENSSEN 2014). Un certain
nombre de ces tests (voir Section \ref{sec:Les-comp=0000E9titeurs-aux tests lisses}) sont disponibles dans le package R \textbf{\noun{BenfordTest}} (JOENSSEN
2013b).

Malgr\'{e} l'importance de l'utilisation d'un \og bon \fg{ }  test d'ad\'{e}quation pour la d\'{e}tection de fraude, ce n'est que r\'{e}cemment (MORROW 2014), (JOENSSEN, 2014), (LESPERANCE et collab. 2016) que les premi\`{e}res analyses de la puissance de ces diff\'{e}rents tests \`{a} la LNB ont \'{e}t\'{e} r\'{e}alis\'{e}es. La proc\'{e}dure pour comparer entre eux plusieurs tests d'ad\'{e}quation est bien rod\'{e}e : on d\'{e}termine une liste d'hypoth\`{e}ses alternatives couvrant les \'{e}carts que l'on estime plausibles \`{a} $H_{0}$, on g\'{e}n\`{e}re des pseudo-\'{e}chantillons de chacune de ces alternatives, on applique les tests au même niveau, soit la probabilit\'{e} d'une erreur de Type I (en g\'{e}n\'{e}ral 0.05), puis on constate s'ils ont rejet\'{e} ou non $H_{0}$.
En r\'{e}p\'{e}tant ceci un grand nombre de fois, on obtient une approximation de la puissance des diff\'{e}rents tests que l'on peut ensuite comparer entre eux. En g\'{e}n\'{e}ral, aucun test n'est uniform\'{e}ment le plus puissant, chacun ayant ses forces et ses faiblesses en regard des alternatives consid\'{e}r\'{e}es. Un  \og bon \fg{ } test se range r\'{e}guli\`{e}rement parmi les tests les plus puissants. Dans ce contexte, l'introduction d'un nouveau test se justifie si on peut montrer qu'il se range aussi parmi les plus puissants pour des alternatives courantes, ou s'il est performant pour de nouvelles alternatives importantes qui n'avaient pas auparavant \'{e}t\'{e} consid\'{e}r\'{e}es.

En d\'{e}tection de fraudes, il s'ajoute \`{a} cette proc\'{e}dure d'\'{e}valuation
de la qualit\'{e} des tests d'ad\'{e}quation une dimension suppl\'{e}mentaire.
Maintenant que l'existence d'outils de d\'{e}tection bas\'{e}s sur la LNB
est bien connue, un fraudeur astucieux va chercher \`{a} trafiquer les
donn\'{e}es de façon \`{a} passer inaperçu (GAUVRIT et collab. 2017). En effet, la crainte d'être d\'{e}tect\'{e} est un puissant frein \`{a} la fraude
(ABDULLAHI et collab. 2015). Si le fraudeur sait que l'auditeur de ses donn\'{e}es trafiqu\'{e}es exploitera tel test d'ad\'{e}quation, il essaiera
de faire en sorte que ce test ne d\'{e}tecte pas d'\'{e}cart \`{a} la LNB. Dans
ce contexte, chaque nouveau test ajoute une contrainte suppl\'{e}mentaire
augmentant la difficult\'{e} de sa t\^{a}che. En outre, avec plusieurs tests
utilis\'{e}s en batterie, l'auditeur a aussi le choix de moduler leur
emploi d'une façon inconnue du fraudeur, augmentant ainsi les risques
de d\'{e}tection. Dans ce contexte, l'introduction d'un nouveau test d'ad\'{e}quation est justifi\'{e}e non seulement par sa bonne puissance, mais aussi parce
que son existence peut complexifier la t\^{a}che du fraudeur. En ce domaine,
viser l'\'{e}radication de la fraude est quasi impossible; on cherche
de façon plus r\'{e}aliste \`{a} la rendre difficile pour aider \`{a} sa
pr\'{e}vention (ABDULLAHI et collab. 2015).

La famille des tests lisses (\emph{smooth tests}) introduite par Neyman (NEYMAN  1937), s'applique \`{a} des donn\'{e}es autant discr\`{e}tes que continues. Ces tests sont plus complexes \`{a} d\'{e}velopper, car ils sont sp\'{e}cifiques \`{a} la loi de probabilit\'{e} postul\'{e}e en $H_{0}$ mais au fil des ann\'{e}es, leurs grandes qualit\'{e}s ont \'{e}t\'{e} reconnues et ceci a men\'{e} (RAYNER et BEST 1990, p. 9) \`{a} la recommandation : \og{Don't use those other methods\textendash use a smooth test}! \fg{ \foreignlanguage{american}{.}
Une version plus moderne, pilot\'{e}e dans les donn\'{e}es, a \'{e}t\'{e} propos\'{e}e (LEDWINA 1994 ) et a conduit (KALLENBERG et LEDWINA 1997 ) \`{a} affiner cette recommandation en :  \og{Use a data-driven smooth test}!\fg{ } . Mais jusqu'\`{a} pr\'{e}sent, ces tests n'ont pas \'{e}t\'{e} d\'{e}velopp\'{e}s pour tester l'ad\'{e}quation \`{a} la LNB. 

Le but du pr\'{e}sent travail est de d\'{e}velopper diff\'{e}rentes variantes
des tests lisses pour le cas de la LNB et d'en \'{e}tudier les qualit\'{e}s.
La Section \ref{sec:Test-lisse-pour la LNB} rappelle les \'{e}l\'{e}ments
th\'{e}oriques permettant de construire une strat\'{e}gie de test lisse et
donne l'expression des statistiques de test pour le cas de la LNB.
Quelques variantes sont introduites ainsi que des r\'{e}sultats th\'{e}oriques
concernant les puissances. La suite du travail explore les avantages
d'inclure ces tests lisses dans une proc\'{e}dure de d\'{e}tection de fraude.
La Section \ref{sec:Les-comp=0000E9titeurs-aux tests lisses} pr\'{e}sente
une liste assez exhaustive des tests d'ad\'{e}quation comp\'{e}titeurs aux
tests lisses et pr\'{e}cise ceux qui sont retenus par la suite. La Section
\ref{sec:Les-alternatives} pr\'{e}sente les alternatives qui sont consid\'{e}r\'{e}es pour la comparaison des puissances des tests retenus. La Section \ref{sec:Simulations} pr\'{e}sente les r\'{e}sultats d'une exp\'{e}rience de simulation qui montre que l'introduction des tests lisses est tout-\`{a}-fait justifi\'{e}e selon les crit\`{e}res \'{e}voqu\'{e}s plus haut.

\section{Test lisse pour la LNB\label{sec:Test-lisse-pour la LNB}}

Le th\'{e}or\`{e}me suivant explique comment construire une famille, index\'{e}e
par l'entier $K$, de tests lisses d'ad\'{e}quation pour l'hypoth\`{e}se nulle
$H_{0}:X\sim f(\cdot)$. Ce th\'{e}or\`{e}me est bien connu et une r\'{e}f\'{e}rence
est (THAS 2010)  par exemple.
\begin{theorem}
\label{theom:1}Soit $X_{1},...,X_{n}$ des copies ind\'{e}pendantes
d'une variable al\'{e}atoire $X$ de densit\'{e} $f(\cdot)$ par rapport \`{a}
une mesure dominante $\nu$. Soit $\{h_{0}(\cdot)\equiv1,h_{k}(\cdot),k=1,2,...\}$
une suite de fonctions orthonormales par rapport \`{a} $f(\cdot)$; plus
pr\'{e}cis\'{e}ment, $\int h_{k}(x)h_{k^{\prime}}(x)f(x)d\nu(x)$ $=\delta_{kk^{\prime}},$
la fonction delta de Kronecker. Soit $U_{k}=n^{-1/2}\sum_{i=1}^{n}h_{k}(X_{i})$
et pour un entier $K\geq1$, soit $T_{K}=\sum_{k=1}^{K}U_{k}^{2}$.\\
Alors sous $H_{0},$ $T_{K}\overset{L}{\longrightarrow}\chi_{K}^{2}$,
la loi khi-deux \`{a} $K$ degr\'{e}s de libert\'{e}, et un test de niveau asymptotique
$\alpha$ rejette $H_{0}$ si la valeur observ\'{e}e de $T_{K}$ d\'{e}passe
$x_{K,1-\alpha}^{2}$, le quantile d'ordre $1-\alpha$ de cette loi
$\chi_{K}^{2}$.
\end{theorem}

Nous sp\'{e}cialisons maintenant ce th\'{e}or\`{e}me au cas où $f(\cdot)$ est
la densit\'{e} de la LNB. Pour ce faire, il faut d\'{e}terminer des fonctions
orthonormales $h_{k}(\cdot)$. Comme tous les moments de la LNB existent,
nous pouvons, \`{a} l'instar de (NEYMAN 1937) et de nombreux autres auteurs par la suite, choisir des polynômes. Le th\'{e}or\`{e}me suivant est aussi
connu  (BOULERICE et DUCHARME 1997). Dans la suite
l'indice \og 0 \fg{ }  d\'{e}note un op\'{e}rateur probabiliste calcul\'{e} sous $H_{0}:X\sim f(\cdot)$.

\begin{theorem}
\label{thm:Construction des polynomes orthogonaux}Soit $\mu_{k}=\mathbb{E}_{0}(X^{k}),k\geq0$.
Soit aussi la matrice $\mathbf{M}_{k}=\left[\mu_{i+i^{\prime}}\right]_{i,i^{\prime}=0,...,k-1}$,
le vecteur $\boldsymbol{\mu}_{k}=(\mu_{k},\mu_{k+1},...\mu_{2k-1})^{T}$
et la constante $c_{k}=\mu_{2k}-\boldsymbol{\mu}_{k}^{T}\mathbf{M}_{k}^{-1}\boldsymbol{\mu}_{k}$. 
Alors les polynômes 
\[ h_{k}(x)=c_{k}^{-1/2}(x^{k}-(1,x,x^{2},...,x^{k-1})\mathbf{M}_{k}^{-1}\boldsymbol{\mu}_{k}) \]
satisfont la condition du Th\'{e}or\`{e}me \ref{theom:1}.
\end{theorem}
Les moments de la LNB ont des expressions explicites complexes. Il
en va de même des coefficients des polynômes $h_{k}(\cdot)$ dont
les expressions exactes sont tr\`{e}s longues d\`{e}s lors que $k>2.$ C'est
pourquoi il est pr\'{e}f\'{e}rable de les exprimer sous une forme approximative.
Mais ces calculs doivent être faits avec soin, car si les approximations
num\'{e}riques sont effectu\'{e}es au niveau des \'{e}l\'{e}ments du Th\'{e}or\`{e}me \ref{thm:Construction des polynomes orthogonaux},
il en d\'{e}coule des erreurs d'arrondis qui d\'{e}truisent l'orthonormalit\'{e}.
Ainsi, les coefficients des $h_{k}(\cdot)$ doivent être calcul\'{e}s
en valeurs exactes, puis convertis en approximations num\'{e}riques. En
utilisant le logiciel \noun{MATHEMATICA}, on obtient ainsi :

\begin{align*}
h_{1}(x) & =-1.3979+0.4063x\\
h_{2}(x) & =2.2836-1.6128x+0.18247x^{2}\\
h_{3}(x) & =4.0815+4.5719x-1.2053x^{2}+0.0862x^{3}\\
h_{4}(x) & =8.0795-12.0946x+5.1951x^{2}-0.8249x^{3}+0.0431x^{4}\\
h_{5}(x) & =-18.1064+33.1385x-19.7207x^{2}+5.0168x^{3}-0.5665x^{4}+0.0233x^{5}
\end{align*}
Le package de R \noun{BenfordSmoothTest}, disponible sur le site \href{https://github.com/st-homme/benfordSmoothTest}{GitHub}, permet de calculer les $T_{K}$ pour $K$ allant jusqu'\`{a} 7. 

\begin{remark}
\label{rem:remarque 3 sur le quantile par Monte Carlo}Dans le Th\'{e}or\`{e}me
\ref{theom:1}, l'approximation $\chi_{K}^{2}$ est bas\'{e}e
sur une convergence quand $n\rightarrow\infty$. Si $n$ est petit,
$x_{K,1-\alpha}^{2}$ peut donner une mauvaise approximation du quantile
exact de la loi de $T_{K}$. S'il est n\'{e}cessaire d'assurer un contrôle
pr\'{e}cis de l'erreur de Type I, on peut approcher la valeur exacte de
ce quantile par la m\'{e}thode de Monte-Carlo. C'est ce qui est fait dans
le Package \noun{BenfordSmoothTest,} d\`{e}s lors que $n$< 100.
\end{remark}

\begin{remark}
\label{rem:remarque 4 sur la version data-driven de Ledwina} Le choix
de l'hyperparam\`{e}tre $K$ est un \'{e}l\'{e}ment important de la strat\'{e}gie
d'un test lisse. Si $K$est trop petit, le test perd de la puissance car il ne peut tenir compte de certains \'{e}carts \`{a} $H_{0}$. Si $K$ est trop grand, cette puissance est dilu\'{e}e par la consid\'{e}ration de termes n\'{e}gligeables dans la statistique de test. Pour choisir judicieusement cet hyperparam\`{e}tre, (LEDWINA 1994) a propos\'{e} une strat\'{e}gie  \og{data - driven} \fg{ } qui consiste \`{a} calculer, pour une valeur \`{a} choisir $K_{max}$,
\begin{ceqn}
\[
\hat{K}=\underset{1\leq k\leq K_{max}}{\arg\max}\left\{ T_{k}-k\:log(n)\right\} ,
\]
\end{ceqn}
et \`{a} utiliser la statistique de test $T_{\hat{K}}$ qui $\overset{L}{\longrightarrow}\chi_{1}^{2}$
sous $H_{0}$. De nombreuses simulations par LEDWINA et ses coauteurs
ont montr\'{e} que ce test \og{data - driven} \fg{ } est un bon compromis dans
la famille de tests lisses $\left\{ T_{k},k=1,...K_{max}\right\} $. Par
ailleurs, on peut choisir $K_{max}$ en exploitant de l'information
contextuelle au probl\`{e}me (ce qu'on appelle le cas \og{horizon fini} \fg{ }),
soit le laisser tendre vers $\infty$ (\og{horizon infini} \fg{}) \`{a} une vitesse
qui d\'{e}pend de $n$. LEDWINA et ses coauteurs ont beaucoup travaill\'{e}
sur cette vitesse et les r\'{e}sultats th\'{e}oriques obtenus sont impressionnants,
mais d'une utilit\'{e} pratique limit\'{e}e, car ils s'expriment sous la forme
$K_{max}=o(n^{j})$ pour un $j$ d\'{e}pendant du contexte du probl\`{e}me.
Heureusement la puissance du test bas\'{e} sur $T_{\hat{K}}$ en fonction
de $K_{max}$ plafonne rapidement, de sorte que le cas  \og{horizon fini} \fg{ }
 avec  $5\leq K_{max}\leq7$ donne en g\'{e}n\'{e}ral de bons r\'{e}sultats.
Nous avons choisi $K_{max}=5$ dans la suite. 
\end{remark}

\begin{remark}
\label{rem:Remqrque sur la puissance du test lisse}Soit $g(x)$ une
alternative fix\'{e}e \`{a} la LNB. La puissance du test bas\'{e} sur $T_{K}$
peut être assez bien approxim\'{e}e par l'expression suivante qui se trouve
dans (\foreignlanguage{american}{INGLOT et collab. 1995},
Th\'{e}or\`{e}me 2.1) et dont les conditions d'applicabilit\'{e} sont rencontr\'{e}es
par virtuellement tous les $g(x)$ raisonnables dans le pr\'{e}sent contexte.
Soit $\nu_{k}=\sqrt{n}\:\sum_{x=1}^{9}h_{k}(x)g(x)$  regroup\'{e}s dans
 $\boldsymbol{\nu}=(\nu_{1},...,\nu_{K})$$^{T}$ et 
$\Sigma=\left[\sigma_{i,j}\right]_{i,j=1,...,K}$,
où $\sigma_{i,j}=\sum_{x=1}^{9}h_{i}(x)h_{j}(x)g(x)-\nu_{i}\nu_{j}/n$.
Posons la d\'{e}composition spectrale $\mathbf{\Sigma}=\mathbf{P}\mathbf{\Lambda}\mathbf{P}^{T}$
où $\mathbf{\Lambda}=\textrm{Diag}\{\lambda_{1},\ldots,\lambda_{K}\}$
est la matrice des valeurs propres de $\Sigma$ et $\mathbf{P}$ et
la matrice de ses vecteurs propres normalis\'{e}s. Alors, uniform\'{e}ment
en $t>0$,
\begin{ceqn}
\begin{equation}
\mathbb{P}_{g}[T_{K}>t]=\mathtt{\mathbb{P}}\left[\sum_{k=1}^{K}\lambda_{k}\chi_{1}^{2}(\delta_{k}^{2})>t\right]+O(n^{-1/2}),
\end{equation}
\end{ceqn}
où $\delta_{k}$ sont les composantes de $\mathrm{\Lambda}^{-1/2}\mathbf{P}\boldsymbol{\nu}$.
Par ailleurs, pour approcher la loi de la somme pond\'{e}r\'{e}e des $\chi_{1}^{2}(\delta_{k}^{2})$,
on peut utiliser une approximation de (LIU et collab. 2009) commod\'{e}ment
bas\'{e}e sur des moments. Adapt\'{e}e au pr\'{e}sent probl\`{e}me, cette approximation
s'\'{e}crit :
\begin{ceqn}
\begin{align*}
\mathtt{\mathbb{P}}\left[\sum_{k=1}^{K}\lambda_{k}\chi_{1}^{2}(\delta_{k}^{2})>x_{K,1-\alpha}^{2}\right] & \simeq\mathtt{\mathbb{P}}\left[\chi_{\ell}^{2}(d^{2})>\frac{a(x_{K,1-\alpha}^{2}-c_{1})}{\sqrt{c_{2}}}+(\ell+d)\right],
\end{align*}
\end{ceqn}
où $\ell=a^{2}-2d^{2}$, $d^{2}=max[s_{1}a^{3}-a^{2},0],$ $a=1/(s_{1}-\sqrt{s_{1}^{2}-s_{2}})$
si $s_{1}^{2}>s_{2}$ et $1/s_{1}$ sinon, avec $s_{1}=c_{3}/c_{2}^{3/2}$
et $s_{2}=c_{4}/c_{2}^{2}$ avec $c_{i}=tr(\Sigma^{i})+i\nu^{T}\Sigma^{(i-1)}\nu$,
$i=1,...,$4. Cette approximation donne en g\'{e}n\'{e}ral de bons r\'{e}sultats
dans le contexte de tests d'ad\'{e}quation (DUCHESNE et LAFAYE DE MICHEAUX
2010). 
L'existence d'une expression explicite approximant la puissance du
test lisse permet de mieux comprendre comment \'{e}volue cette puissance
en fonction de $g(x)$. En effet, pour plusieurs des tests concernant
la LNB (voir Section \ref{sec:Les-comp=0000E9titeurs-aux tests lisses}),
la seule façon de calculer la puissance est par simulations, ce qui
ne permet pas une bonne \'{e}tude de la sensibilit\'{e} de cette puissance
en regard des \'{e}carts entre $f(x)$ et $g(x)$.
\end{remark}

\begin{remark}
\label{rem:Remarque sur la stat a-star} Dans la famille des tests
lisses $\left\{ T_{k},k=1,...K_{max}\right\} $, on retrouve souvent
des statistiques d\'{e}j\`{a} propos\'{e}es dans la litt\'{e}rature sur la base de
consid\'{e}rations intuitives. C'est le cas ici, où $\left|h_{1}(x)\right|$
/2.25915 = $a^{*}$, la statistique introduite par (JUDGE et SCHECHTER 2009) et que l'on retrouve dans le package R \noun{BenfordTest }(fonction
\noun{meandigit.benftest}). Cette statistique est aussi reli\'{e}e \`{a} la
statistique du \og Distortion Factor (DF) model \fg{ } de (NIGRINI 1996) .
\end{remark}

\section{Les comp\'{e}titeurs aux tests lisses\label{sec:Les-comp=0000E9titeurs-aux tests lisses}}
Comme signal\'{e} dans l'introduction, le nombre de tests d'ad\'{e}quation \`{a} la LNB est plutôt limit\'{e} au del\`{a} du test du $\chi^{2}$ de Pearson de (\ref{eq:Pearson chi-square test}) et ceci facilite la t\^{a}che des fraudeurs. Soit $D_{1},...,D_{n}$ un \'{e}chantillon al\'{e}atoire de \emph{PCS}. Pour $d=1,..,9$, d\'{e}notons par $\hat{p}_{d}=n_{d}/n$ , la proportion de $d$ dans l'\'{e}chantillon et $\pi_{d}$ les probabilit\'{e}s de la LNB. 
Soit aussi $S_{d}=\sum_{j=1}^{d}\hat{p}_{j}$ et $S_{d}^{*}=\sum_{j=1}^{d}\pi_{j}$.
Posons $Z_{d}=S_{d}-S_{d}^{*}$ et $t_{d}=(\pi_{d}+\pi_{d+1})/2$,
pour $d=1,..,8$ et $t_{9}=(\pi_{9}+\pi_{1})/2$. (LESPERANCE et collab. 2016) consid\`{e}rent les versions discr\`{e}tes des tests suivants bas\'{e}s
sur les \'{e}carts entre la distribution cumulative empirique et la distribution
cumulative th\'{e}orique de la LNB. 
\begin{ceqn}
\begin{align*}
W_{n}^{2} & =n\sum_{i=1}^{n}Z_{d}^{2}t_{d}\qquad\textrm{(Cramer-von Mises)},
\end{align*}
\end{ceqn}
\begin{ceqn}
\begin{align*}
U_{n}^{2} & =n\sum_{i=1}^{n}(Z_{d}-\bar{Z})^{2}t_{d},\qquad\textrm{(Watson)},
\end{align*}
\end{ceqn}
dont une variante, dite de Freedman, se retrouve dans (JOENSSEN 2014)
et
\begin{ceqn}
\begin{align*}
A^{2} & =n\sum_{i=1}^{8}\frac{Z_{d}^{2}t_{d}}{T_{d}(1-T_{d})}\qquad\textrm{(Anderson-Darling)}.
\end{align*}
\end{ceqn}
(MORROW 2014) consid\`{e}re la statistique de Kolmogorov $K_{n}=\sqrt{n}\underset{1\leq d\leq9}{\max}\left|S_{d}-S_{d}^{*}\right|.$
Enfin, un certain nombre d'auteurs consid\`{e}rent des tests bas\'{e}s sur
les \'{e}carts entre les probabilit\'{e}s $\pi_{d}$ de la LNB et les $\hat{p}_{d}$.
(LEEMIS et collab. 2000) proposent la statistique $m=\underset{1\leq d\leq9}{\max}\left|\hat{p}_{d}-\pi_{d}\right|$
alors que (CHO et GAINES 2007) sugg\`{e}rent $d=\sqrt{\sum_{d=1}^{9}(\hat{p}_{d}-\pi_{d})^{2}}$ et (DRAKE et NIGRINI 2000) introduisent la \og Mean Average Deviation \fg{ } : $MAD=\frac{1}{9}\sum_{d=1}^{9}\left|\hat{p}_{d}-\pi_{d}\right|$.
Enfin, dans un autre ordre d'id\'{e}e, (JUDGE et SCHECHTER 2009) proposent la statistique $a^{*}=(\bar{D}-3.44027)/(5.55973)$ qui est proportionnelle
\`{a} la statistique $T_{1}$ du test lisse (voir Remarque \ref{rem:Remarque sur la stat a-star}).
Au meilleur de nos connaissances, cette nomenclature couvre pratiquement tous les tests existants \`{a} ce jour pour la LNB, sauf la statistique $J_{p}^{2}$ de (JOENSSEN 2013a) inspir\'{e}e du test de Shapiro-Wilks, qui ne sera pas consid\'{e}r\'{e} plus avant en raison de son mauvais comportement, un test adaptant la statistique de Hotelling qui semble peu int\'{e}ressant car certains hyperparam\`{e}tres doivent être choisis (sans crit\`{e}re pr\'{e}cis \`{a} ce jour), et un test bay\'{e}sien  (GEYER et WILLIAMSON 2004) qui ne permet pas le contr\^{o}le fr\'{e}quentiste des erreurs de Type I et II, \'{e}l\'{e}ment important pour les auditeurs. 

Le package R \noun{BenfordTests} (version 1.2.0, 2015) maintenu par
(JOENSSEN 2013b) permet d'effectuer certains de ces tests via les
fonctions\noun{ chisq.benftest} (test du $\chi^{2}$),\noun{ ks.benftest}
(test de Kolmogorov), \noun{mdist.benftest }(test bas\'{e} sur \emph{m})
, \noun{edist.benftest} (test bas\'{e} sur \emph{d}), \noun{usq.benftest}
(variante de Freedman du test de Watson), \noun{jpsq.benftest} (test
bas\'{e} sur $J_{p}^{2}$) et \noun{meandigit.benftest} (test bas\'{e} sur
$a^{*}$), auquel s'ajoute le test adaptant la statistique 
de Hotelling.

Signalons que (LESPERANCE et collab. 2016)
et (JOENSSEN 2014) recommandent le test de Watson ou sa variante de
Freedman qui donnent de bons r\'{e}sultats  dans leurs simulations dont
nous reprenons certains \'{e}l\'{e}ments aux Sections \ref{sec:Les-alternatives}
et \ref{sec:Simulations}, en ceci que leur puissance se range de
façon consistante parmi les plus \'{e}lev\'{e}es pour les quelques alternatives
qu'ils consid\`{e}rent. Dans la suite de ce travail, nous comparons la
puissance des pr\'{e}sents tests lisses \`{a} celles de certains des tests
plus haut sur une plage plus large d'alternatives. Plus pr\'{e}cis\'{e}ment, comme repr\'{e}sentant de la famille des tests lisses $\left\{ T_{k},k=1,...K_{max}\right\} $, nous retenons le test $T_{2}$ et sa version data-driven $T_{\hat{K}}$ (avec $K_{max}=5)$; comme repr\'{e}sentant des tests bas\'{e}s sur les \'{e}carts entre distributions cumulatives, nous choisissons la statistique $U_{n}^{2}$  de Watson. Nous consid\'{e}rons le test MAD de (DRAKE et NIGRINI 2000) comme repr\'{e}sentant des tests bas\'{e}s sur des \'{e}carts entre les $\pi_{d}$ de la LNB et les $\hat{p}_{d}$, auquel nous ajoutons le  test classique du $\chi^{2}$ de Pearson (\ref{eq:Pearson chi-square test}). 

\section{Les alternatives\label{sec:Les-alternatives}}

Pour \'{e}tudier la puissance des tests lisses de la Section \ref{sec:Test-lisse-pour la LNB} et la comparer avec celle de quelques tests repr\'{e}sentatifs de la Section \ref{sec:Les-comp=0000E9titeurs-aux tests lisses}, nous devons pr\'{e}ciser des alternatives pour lesquelles cette puissance sera calcul\'{e}e. Pour ce faire, nous consid\'{e}rons des familles d'alternatives index\'{e}es par un param\`{e}tre, g\'{e}n\'{e}riquement not\'{e} $\beta$, et emboîtant la LNB. Notre premi\`{e}re famille d'alternatives est celle de (RODRIGUEZ 2004) donn\'{e}e par :

\begin{ceqn}
\begin{align}
\mathbb{P}[D=d]=p_{d}^{(Rod)}(\beta) & =\begin{cases}
\frac{1}{9}(1+\frac{10}{9}\ln(10)+x\ln(x)-(x+1)\ln(x+1)) & \textrm{si }\beta=0\\
\log_{10}[1+1/x] & \textrm{si }\beta=-1\\
\frac{\beta+1}{9\beta}-((x+1)^{(\beta+1)}-x^{(\beta+1)})/(\beta(10^{(\beta+1)}-1)) & \textrm{sinon}
\end{cases},\label{eq:Rodriguez family}
\end{align}
\end{ceqn}
pour $\beta\in\mathbb{R}$, où on peut remarquer que quand $\beta=0$,
on retrouve la loi de Stigler (LEE et collab. 2010), $\beta=-1$
donne la LNB et quand $\beta\rightarrow\pm\infty$, on a la loi uniforme
discr\`{e}te sur $\{1,2,...,9\}$, not\'{e}e $UD[\{1,..,9\}]$. La deuxi\`{e}me
famille d'alternatives est celle de la LNB g\'{e}n\'{e}ralis\'{e}e de (PIETRONERO 
et collab. 2001) donn\'{e}e par : 

\begin{ceqn}
\begin{align*}
\mathbb{P}[D=d] & =p_{d}^{(Piet)}(\beta)=\begin{cases}
\frac{(d+1)^{(1-\beta)}-d{}^{(1-\beta)}}{10^{(1-\beta)}-1} & \textrm{si }\beta\neq1\\
\log_{10}[1+1/x] & \textrm{si }\beta=1
\end{cases},
\end{align*}
\end{ceqn}
où $\beta\in\mathbb{R}$ avec la LNB correspondant au cas $\beta=1$.
La troisi\`{e}me famille de lois alternatives est celle de (HURLIMANN 2006)
où

\begin{alignat*}{1}
\mathbb{P}[D=d]=p_{d}^{(Hurl)}(\beta)= & \frac{1}{2}(\log_{10}[1+x]^{\beta}-\log_{10}[x]^{\beta}-(1-\log_{10}[1+x])^{\beta}+(1-\log_{10}[x])^{\beta}),
\end{alignat*}
avec $\beta\in\mathbb{R}$. Notons que cette famille a la particularit\'{e} de donner la LNB pour les deux valeurs $\beta=1,2$. La quatri\`{e}me famille est celle d'un m\'{e}lange de lois LNB et UD de la forme $(1-\beta)\times LNB+\beta\times UD[\{0,1,...,9\}]$, avec $\beta\in[0,1],$ consid\'{e}r\'{e}e par (LESPERANCE et collab. 2016). La cinqui\`{e}me famille de lois, aussi consid\'{e}r\'{e}e par (LESPERANCE et collab. 2016), est celle d'une loi LNB contamin\'{e}e où 

\begin{align*}
\mathbb{P}[D=d] & =p_{d}^{(conta-1)}(\beta)=\begin{cases}
\pi_{d}/(1+2\beta) & \textrm{si }d\neq1,9\\
(\pi_{d}+\beta)/(1+2\beta) & \textrm{si }d=1,9
\end{cases},
\end{align*}
avec $\beta\in[0,0.6]$. Enfin la derni\`{e}re famille de lois est nouvelle;
elle est obtenue en contaminant plus que la famille pr\'{e}c\'{e}dente l'aile
de droite de la LNB de la façon suivante : 

\begin{align*}
\mathbb{P}[D=d] & =p_{d}^{(conta-2)}(\beta)=\begin{cases}
\pi_{d}/(1+10\beta) & \textrm{si }d=1,...,5\\
(\pi_{d}+(d-5)\beta)/(1+10\beta) & \textrm{si }d=6,...,9
\end{cases},
\end{align*}
où $\beta=0.001(1+\gamma/2)$, avec la LNB correspondant au cas $\gamma=-2$.
Ensemble, ces 6 familles couvrent de nombreuses alternatives \`{a} la
LNB. Elles ont aussi \'{e}t\'{e}, pour les cinq premi\`{e}res, reconnues comme \'{e}tant des
alternatives plausibles, pouvant être rencontr\'{e}es dans les applications.
L'\'{e}tude de la puissance des diff\'{e}rents tests d'ad\'{e}quation devrait
donner une bonne id\'{e}e de leur comportement en pratique. Terminons
cette section en signalant que les cinq premi\`{e}res familles sont, \`{a}
notre connaissance, les seules familles de lois existantes dans la litt\'{e}rature qui emboîtent
la LNB.

\section{Simulations\label{sec:Simulations} }

Dans cette section, nous \'{e}tudions la puissance de diff\'{e}rents tests
pour les alternatives de la Section \ref{sec:Les-alternatives}. Pour chacune des 6 familles d'alternatives de la Section \ref{sec:Les-alternatives}, nous avons choisi une taille d'\'{e}chantillon repr\'{e}sentative des comportements de la fonction de puissance, tout en \'{e}vitant les r\'{e}sultats triviaux (puissance trop proche de $\alpha=0.05$ et 1). Ainsi pour la famille de Rodriguez, $n=250$; pour celle de Pietronero, $n=50$; pour la famille de Hurlimann, $n=750$ et enfin pour la famille de m\'{e}lange de LNB et UD ainsi que les deux LNB contamin\'{e}es, nous prenons $n=500$. Les tests sont effectu\'{e}s au niveau 5\% et en regard de la Remarque \ref{rem:remarque 3 sur le quantile par Monte Carlo}, les quantiles de r\'{e}f\'{e}rence de tous les tests sont approxim\'{e}s par Monte Carlo (50 000 r\'{e}plications) afin de permettre une juste comparaison des puissances. Ensuite, pour chacune des familles, et en faisant varier la valeur du param\`{e}tre $\beta$ sur une plage de valeurs adapt\'{e}es
\`{a} chaque famille, 10 000 pseudo-\'{e}chantillons des tailles prescrites
sont g\'{e}n\'{e}r\'{e}s. Les tests sont appliqu\'{e}s \`{a} chacun des pseudo-\'{e}chantillons et les puissances de chaque test sont approxim\'{e}es par le nombre de rejets parmi ces 10 000 r\'{e}plications. Notons que pour le test $T_{2}$, les formules de la Remarque \ref{rem:Remqrque sur la puissance du test lisse}
sont aussi calcul\'{e}es et produisent des r\'{e}sultats remarquablement proches
de ceux de la simulation. 

Les r\'{e}sultats de cette exp\'{e}rience apparaissent \`{a} la Figure \ref{fig:Puissance-de-divers}.
Comme les puissances de cette figure sont approxim\'{e}es de 10 000 r\'{e}plications,
elles contiennent un bruit statistique que l'on peut \'{e}valuer \`{a} environ
$\pm0.01$ (au niveau de confiance 95\%) lorsque la puissance est
autour de 0.5. Ce point est pris en consid\'{e}ration dans les conclusions
qui suivent.

Pour les lois de Rodriguez du Panneau a), le test de Watson bas\'{e} sur
$U_{n}^{2}$ est l\'{e}g\`{e}rement sup\'{e}rieur \`{a} $T_{2}$ quand $\beta<-1$,
qui \`{a} son tour est suivi de pr\`{e}s par $T_{\hat{K}}$. Quand $\beta>-1$,
ces trois tests ont approximativement la même puissance. Les tests
bas\'{e}s sur MAD et $\chi^{2}$ sont moins puissants (avec MAD meilleur
que $\chi^{2}$ si $\beta<-1$ et inversement quand $\beta>-1$).
Pour les lois de Pietronero du Panneau b), $T_{2}$ est le test le
plus puissant sur toute la plage de valeurs de $\beta$, suivi de
$T_{\hat{K}}$ puis de $U_{n}^{2}$ quand $\beta<1$ et inversement
quand $\beta>1$. Encore l\`{a}, les tests MAD et $\chi^{2}$ ferment
la marche, comparables quand $\beta<1$ mais avec $\chi^{2}$ nettement
moins puissant que MAD \`{a} droite. Pour la famille de Hurlimann du Panneau
c), quand $\beta<1,T_{2}$ domine, suivi de $T_{\hat{K}}$ puis
de $U_{n}^{2}$ alors que pour $\beta>1,$ $T_{2}$ et $U_{n}^{2}$
sont comparables mais $T_{\hat{K}}$ est l\'{e}g\`{e}rement moins puissant.
Encore l\`{a}, MAD et $\chi^{2}$ ferment la marche. Pour la famille de
m\'{e}lange LNB et UD du Panneau d), $T_{\hat{K}}$ domine $T_{2}$, MAD
et $\chi^{2}$ sont les moins puissants et $U_{n}^{2}$ a un comportement
interm\'{e}diaire. Pour la premi\`{e}re famille de lois contamin\'{e}es du Panneau e), $T_{2}$ et $\chi^{2}$ dominent les autres tests, suivi de pr\`{e}s
par $T_{\hat{K}}$ alors que $U_{n}^{2}$ et MAD ont une puissance
nettement plus faible. Enfin pour la famille des lois contamin\'{e}es
du Panneau f), les deux tests lisses $T_{2}$ et $T_{\hat{K}}$ dominent
nettement les autres tests, MAD et $U_{n}^{2}$ \'{e}tant de loin inf\'{e}rieurs
alors que $\chi^{2}$ a un comportement interm\'{e}diaire. 

En conclusion, on retire de cette exp\'{e}rience de simulation que le
test $T_{2}$ est toujours parmi les deux meilleurs tests en terme
de puissance, disputant souvent la tête avec  le test bas\'{e} sur $T_{\hat{K}}$.
Le test MAD est toujours parmi les deux plus mauvais et ne devrait
pas être utilis\'{e}. La puissance des tests du $\chi^{2}$ et $U_{n}^{2}$
a un comportement plus variable. En regard de notre d\'{e}finition de
 \og{bon}\fg{ } test d'ad\'{e}quation \'{e}voqu\'{e}e dans l'introduction, nous pouvons
donc recommander l'un ou l'autre des tests $T_{2}$ ou $T_{\hat{K}}$.
Leur introduction est ainsi parfaitement justifi\'{e}e et le fait qu'il
s'agisse de tests nouveaux, augmentant l'arsenal des auditeurs, ajoute
\`{a} leur utilit\'{e}. Nous confirmons (MORROW 2014) que le test
du $\chi^{2}$ de Pearson est g\'{e}n\'{e}ralement inf\'{e}rieur et rangeons
le test MAD dans la même cat\'{e}gorie. Enfin, nos r\'{e}sultats att\'{e}nuent
les conclusions de (LESPERANCE et collab. 2016) et (JOENSSEN 2014) qui recommandent 
le test $U_{n}^{2}$ de Watson.

\begin{figure}[htbp]
\centerline{\includegraphics[scale=0.8]{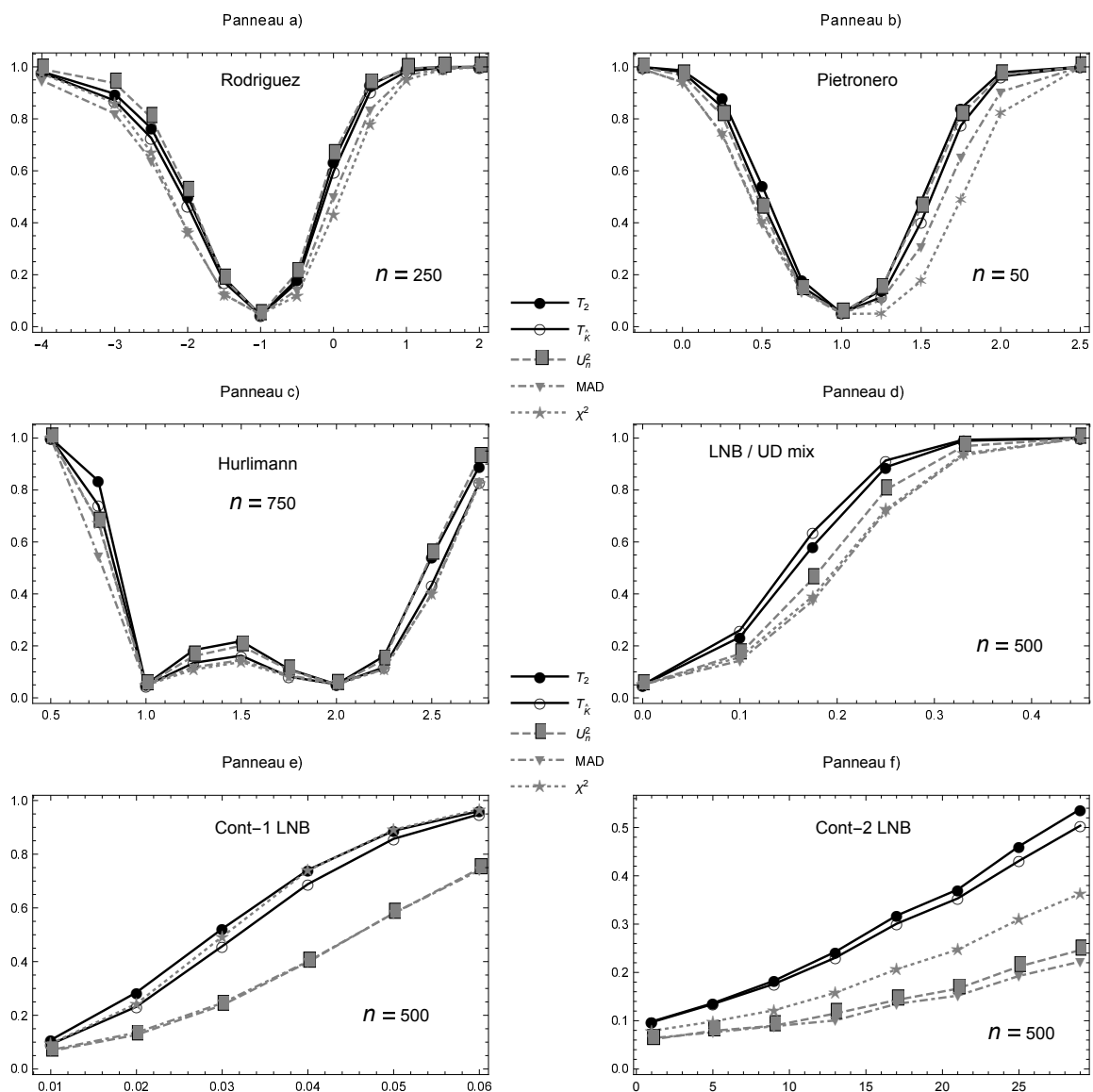}}\caption{\textit{\sffamily Puissance de divers tests (bas\'{e}s sur 10 000 r\'{e}plications) au niveau 5\% pour l'hypoth\`{e}se nulle de loi LNB.
Les familles d'alternatives sont d\'{e}crites \`{a} la Section \ref{sec:Les-alternatives}.
Les tests consid\'{e}r\'{e}s sont $T_{2}$, (trait plein avec cercle plein), $T_{\hat{K}}$ (trait plein avec cercle vide), $U_{n}^{2}$ (tiret avec rectangle), MAD (tiret-pointill\'{e} avec triangle) et le test $\chi^{2}$ de Pearson (pointill\'{e} avec \'{e}toile), dont les expressions se trouvent \`{a} la Section \ref{sec:Les-comp=0000E9titeurs-aux tests lisses}. Les quantiles de r\'{e}f\'{e}rence sont approxim\'{e}s par Monte Carlo en utilisant 50 000 r\'{e}p\'{e}titions.}}
\label{fig:Puissance-de-divers}
\end{figure}

\end{document}